%% file: main.tex
\documentclass[
oneside,
headsepline,
10pt,
BCOR=10mm,
titlepage=false]{scrartcl}

\usepackage{amsfonts}
\usepackage{amssymb}
\usepackage{amsthm}
\usepackage{amsmath, amsfonts, amsxtra}
\usepackage{verbatim} %
\usepackage{esint}
\input{environments}

\input{commands}

\numberwithin{equation}{section}%

\title{A Note on Regularity for the $n$-dimensional $H$-System assuming logarithmic higher Integrability}

\author{Armin Schikorra}

\usepackage[bookmarks=true,linktocpage=true,pdftex=true,pdfpagelabels=true,plainpages=false,pagebackref=false]{hyperref}	

\hypersetup{
pdfauthor={Armin Schikorra},
pdftitle={DAMB},
pdfstartview={FitV},
pdfstartpage={1},
pdfhighlight={/I},
pdfborder={0 0 0 0},
bookmarksopen=false,
bookmarksnumbered=true,
pdfpagemode={UseOutlines},
hypertexnames=true,
unicode=false,
breaklinks=true,
frenchlinks=false,
colorlinks=true,
linkcolor=blue,
citecolor=blue,
urlcolor=blue,
filecolor=blue,
pdftoolbar=true,
pdfmenubar=true,
pdfcenterwindow=true,
pdffitwindow=true,
pdfnewwindow=true
}

\begin{document}
\maketitle
\begin{abstract}
\noindent 
We prove H\"older continuity for solutions to the $n$-dimensional $H$-System assuming logarithmic higher integrability of the solution.
\\[1ex]
{\bf Keywords:} Harmonic maps, nonlinear elliptic PDE, regularity of solutions.\\
{\bf AMS Classification:} 58E20, 35B65, 35J60, 35S05.  	
\end{abstract}
\tableofcontents
\thispagestyle{empty}
\input{intro}
\newpage
\input{orlicz}
\newpage
\input{frehse}
\newpage
\input{pdes}
\newpage

\bibliographystyle{alpha}%
\bibliography{bib}%
\vspace{7em}
\begin{tabbing}
\quad\=Armin Schikorra\\
\>Institut f\"ur Mathematik\\
\>RWTH Aachen University\\
\>Templergraben 55\\
\>52062 Aachen\\
\>Germany\\
\\
\>email: schikorra@instmath.rwth-aachen.de\\
\>page: www.instmath.rwth-aachen.de/$\sim$schikorra
\end{tabbing}
\end{document}

%% file: environments.tex
\newtheorem{theorem}{Theorem}[section]%
\newtheorem{lemma}[theorem]{Lemma}%
\newtheorem{corollary}[theorem]{Corollary}%



\makeatletter

\newenvironment{proofT}[1]{\par
  \pushQED{{\small\it Theorem {\rm #1}} {\leavevmode
  \hbox to.77778em{%
  \hfil\vrule
  \vbox to.675em{\hrule width.6em\vfil\hrule}%
  \vrule\hfil}}}%
  \normalfont \topsep6\p@\@plus6\p@\relax
  {\it \textbf{\proofname { of Theorem {\rm #1}\@addpunct{.}}}\newline\rm}
}{%
   \begin{flushright}\popQED\end{flushright}\@endpefalse
}
\makeatother

\makeatletter

\newenvironment{proofL}[1]{\par
  \pushQED{{\small\it Lemma {\rm #1}} {\leavevmode
  \hbox to.77778em{%
  \hfil\vrule
  \vbox to.675em{\hrule width.6em\vfil\hrule}%
  \vrule\hfil}}}%
  \normalfont \topsep6\p@\@plus6\p@\relax
  {\it \textbf{\proofname { of Lemma {\rm #1}\@addpunct{.}}}\newline\rm}
}{%
   \begin{flushright}\popQED\end{flushright}\@endpefalse
}
\makeatother

\makeatletter

\newenvironment{proofC}[1]{\par
  \pushQED{{\small\it Corollary {\rm #1}} {\leavevmode
  \hbox to.77778em{%
  \hfil\vrule
  \vbox to.675em{\hrule width.6em\vfil\hrule}%
  \vrule\hfil}}}%
  \normalfont \topsep6\p@\@plus6\p@\relax
  {\it \textbf{\proofname { of Corollary {\rm #1}\@addpunct{.}}}\newline\rm}
}{%
   \begin{flushright}\popQED\end{flushright}\@endpefalse
}
\makeatother

\usepackage{array} 
\newenvironment{ma}{\begin{array}{>{\displaystyle}r >{\displaystyle}c >{\displaystyle}l}}{\end{array}}%

%% file: commands.tex
\newcommand{\aleq}{\prec}
\newcommand{\ageq}{\succ}
\newcommand{\aeq}{\approx}

\newcommand{\R}{{\mathbb R}}

\newcommand{\fracm}[1]{\frac{1}{#1}}
\newcommand{\dv}{\operatorname{div}}
\newcommand{\curl}{\operatorname{curl}}
\newcommand{\dn}[1]{\dv(\abs{#1}^{n-2} #1)}

\newcommand{\EXP}{\operatorname{EXP}}%

\newcommand{\lap}{\Delta}

\newcommand{\abs}[1]{{\left \vert #1 \right \vert}}%
\newcommand{\brac}[1]{{\left ( #1 \right )}}%

\newcommand{\intl}{\int \limits}

\newcommand{\mvint}{\fint}

\def\XXint#1#2#3{{\setbox0=\hbox{$#1{#2#3}{\int}$}%
     \vcenter{\hbox{$#2#3$}}\kern-.5\wd0}}%

\newcommand{\sref}[2]{#1.\ref{#2}}


%% file: intro.tex
\section{Introduction}
In his seminal work \cite{Riv07} Tristan Rivi\`{e}re proved that for domains $D \subset \R^2$ solutions $u \in W^{1,2}(D,\R^N)$ to the equation
\begin{equation}\label{eq:rivsys}
 \dv(\nabla u_i) = \Omega_{ik} \cdot \nabla u_k \qquad \mbox{in $D$}
\end{equation}
are H\"older-continuous if $\Omega_{ik} \in L^2(D)$ and $\Omega_{ik} = - \Omega_{ki}$, $1 \leq i,k \leq N$. This gain in regularity (as the right-hand side of the equation belongs a priori to $L^1(D)$, no classic approach implies better regularity of $u$) relies heavily on the antisymmetry of $\Omega$ and is closely related to the relation between Hardy-spaces and div-curl-quantities \cite{CLMS93}. The importance of this result is that equation \eqref{eq:rivsys} is a model equation for Euler-Lagrange equations of critical conformally invariant variational functionals in two dimensions. We refer to the exhaustive introduction in \cite{Riv07} for details.\\
A possible extension of this result, as was suggested in \cite{Riv08}, is the case of the $n$-Laplacian ($n \geq 2$), that is $u \in W^{1,n}(D,\R^N)$ a solution to
\begin{equation}\label{eq:rivsysnD}
 \dv(\abs{\nabla u}^{n-2}\nabla u_i) = \abs{\nabla u}^{n-2}\ \Omega_{ik} \cdot \nabla u_k \quad \mbox{in $D \subset \R^n$}.
\end{equation}
Note that this equation reduces to \eqref{eq:rivsys} if $n=2$. It is also a possible model for a variety of geometrically motivated equations such as the $n$-harmonic maps and the $H$-system. There are some regularity results if additional integrability conditions on the right-hand side are imposed, see e.g. \cite{IO07}, \cite{DM10}, yet these do not use a compensation phenomenon as div-curl-products or antisymmetry, i.e. the assumptions are rather strong. Other results, e.g., \cite{Strz94}, \cite{Fuchs93}, treat the special case of $n$-harmonic maps into spheres. The starting point of this note is yet another kind of result by Kolasi\'{n}ski \cite{Kol09} in the special case of the $n$-dimensional $H$-System. It relies on an additional condition of the differentiability of the solution, but also uses crucially a div-curl-compensation effect:\\
Assume $H \in W^{1,\infty}(\R^{n+1})$ and let $u \in W^{1,n}(D,\R^{n+1})$ be a solution to
\begin{equation}\label{eq:Hsys}
 \dv(\abs{\nabla u}^{n-2}\nabla u_i) = H(u)\ \brac{u_{x_1} \times \ldots \times u_{x_n}}_i \quad \mbox{in $D \subset \R^n$}.
\end{equation}
Here $\times$ is the usual cross product for vectors in $\R^{n+1}$. If one assumes that $u \in W^{n-1,n'}(D,\R^{n+1})$ then $u$ is H\"older continuous.\\
Technically, the proof in \cite{Kol09} relies on growth estimates of local $L^p$-norms, $p < n$, of the gradient of the solutions on small balls, in order to apply Dirichlet Growth Theorem. 
\\
In this note we are concerned with replacing the additional \emph{differentiability} condition $u \in W^{n-1,n'}$ by an additional \emph{integrability} condition of $\nabla u$, and an intuitive approach towards that kind of result might look like this:\\
Note that $u \in W^{n-1,n'}(D,\R^{n+1})$ implies in particular that $\nabla u \in L^{n,n'}$. The latter space is a Lorentz-space which is a strict subspace of $L^n$, cf. \cite{Hunt66}, \cite{GrafC08}.\\
If one observes how the estimates in \cite{Kol09} of the $L^p$-Norm on small balls of the gradient of a solution to \eqref{eq:Hsys} behave if $p$ tends to $n$, one (na\"ively) might be tempted to conjecture that a sufficient condition for regularity might be the integrability $\nabla u \in L^{n,n'}$ or even better $\nabla u \in L^{n,2}$: one estimates the growth of the $L^{n,\infty}$-Norm of the gradient of the solution (assuming that this implies estimating the right-hand side of \eqref{eq:Hsys} tested by functions $\varphi$ with $\nabla \varphi$ bounded in $L^{n,1}$).
Note that by this kind of argument -- if it worked -- one could also obtain a similar result as in \cite{DM10}, where \emph{no structure} but a rather strict \emph{higher integrability} on the right hand side is assumed.\\
But as it turns out, due to the nonlinearity of the $n$-Laplace, the growth of the $L^{n,\infty}$-Norm of $\nabla u$ on small sets does not seem to be that easily estimated by the $n$-Laplace of $u$.\\
Another possibility to replace the differentiability condition $u \in W^{n-1,n'}$ by an integrability-condition is to use logarithmic Orlicz spaces (for the relevant definitions see Section~\ref{sec:orlicz}). If one assumes that $\nabla u \in L^n \log^\alpha L$, $\alpha \in [0,n-1)$ then still no standard growth condition implies that a solution to \eqref{eq:Hsys} is continuous. Nevertheless,
\begin{theorem}\label{th:th}
There exists $\varepsilon > 0$ such that the following holds. Assume $D \subset \R^n$ and let $u \in W^{1,n}(D,\R^{n+1})$ be a solution to the following equation in $D$
\begin{equation}\label{eq:hmodel}
 \dn{\nabla u_i} = H\ \sum_{k=1}^{n+1} \lambda_{k,i}\ \det (\nabla u_1, \ldots, \nabla u_{k-1},\nabla u_{k+1}, \ldots, \nabla u_{n+1}).
\end{equation}
Here, $\lambda_{k,i} \in \R$, $1 \leq i,k \leq n+1$. If $H \in L^\infty \cap W^{1,n}(\R^{n+1},\R)$, and \emph{moreover} $\nabla H, \nabla u \in L^n \log^{n-1-\varepsilon} L$, then $u \in C^{0,\beta}$ for some $\beta > 0$.\\
\end{theorem}
Note that if $\nabla u \in L^n\log^{n-1-\varepsilon}L$ then the right-hand side of \eqref{eq:hmodel} belongs to $L^1 \log^{n-1-\varepsilon} L$, i.e. we are below the range of continuity in \cite{IO07}, where the right-hand side belongs to $L^1 \log^{n-1+\varepsilon} L$. In fact, Frehse's counterexample still holds (see Section~\ref{s:frehse}), i.e. $\nabla \log \log \frac{4}{\abs{x}} \in L^n \log^{n-1-\varepsilon} L$ for any $\varepsilon > 0$.\\
In particular, this theorem implies (see, e.g., \cite{Kol09}) that solutions to the $n$ dimensional $H$-System \eqref{eq:Hsys} are H\"older continuous, if one assumes that $\nabla u \in L^n \log^{n+1-\varepsilon} L$. We would like to stress that Kolasi\'{n}ski's result in \cite{Kol09} and our's  are complementary: Neither contains the other as a special case.\\
\\
Our proof relies on the observation, that with our integrability assumptions it is possible to test \eqref{eq:hmodel} with (a suitably mollified version of) the solution itself. Indeed, testing \eqref{eq:hmodel} with $\varphi \in C_0^\infty(D)$, we have (roughly) by the Hardy-BMO-inequality and the results of \cite[Theorem II.1]{CLMS93}
\begin{equation}\label{eq:rhsest}
\begin{ma}
&&\int_D \abs{\nabla u}^{n-2} \nabla u\ \cdot \nabla \varphi\\
&\leq& C\ \brac{\Vert H \Vert_\infty\ \Vert \nabla \varphi \Vert_{n} + \Vert \varphi \nabla H \Vert_n}\ \Vert \nabla u \Vert_{n}^n.
\end{ma}
\end{equation}
The only term posing a potential problem if $\varphi \to u$ is
\[
 \Vert \varphi \nabla H \Vert_n.
\]
But there is a well-known duality between $L\log L$ and $\EXP$ (see Lemma \ref{la:dualexplnlogalv2} in Section~\ref{sec:orlicz}, or \cite{BS88}), and one has
\[
 \Vert \varphi \nabla H \Vert_n^n \leq \Vert \abs{\varphi}^n \Vert_{\EXP}\ \Vert \abs{\nabla H}^n \Vert_{L \log L} =
\Vert \varphi \Vert_{\EXP,n}^n\ \Vert \nabla H \Vert_{n, \log 1}^n.
\]
Of course, similar inequalities hold for the space $L \log^{n-1} L$, and we can compute
\[
 \Vert \varphi \nabla H \Vert_n^n \leq C\ \Vert \varphi \Vert_{\EXP,\frac{n}{n-1}}^n\ \Vert \nabla H \Vert_{n, \log(n-1)}^n.
\]
Finally we have Trudinger's inequality,
\[
 \Vert \varphi \Vert_{\EXP,\frac{n}{n-1}} \leq C \Vert \nabla \varphi \Vert_{n}.
\]
Together this implies,
\[
 \int_D \abs{\nabla u}^{n-2} \nabla u\ \cdot \nabla \varphi \leq C \brac{\Vert H \Vert_\infty,\Vert \nabla H \Vert_{n,\log(n-1)}}\ \Vert \nabla \varphi \Vert_{n}\ \Vert \nabla u \Vert_n^n.
\]
Localizing this argument to small balls and using Dirichlet Growth Theorem, one concludes regularity for $u$ if $\nabla H \in L^n \log^{n-1} L$. In order to relax the condition
\[
 \Vert \nabla H \Vert_{n,\log(n-1)} < \infty
\]
into 
\[
 \Vert \nabla H \Vert_{n,\log(n-1-\varepsilon)} < \infty, \qquad \mbox{(for some $\varepsilon > 0$)}
\]
we use an adaption of the result of \cite[Theorem II.1]{CLMS93} in Orlicz-Spaces, to have instead of \eqref{eq:rhsest} for some $\alpha, \beta > 0$
\[
\begin{ma}
&&\int_D \abs{\nabla u}^{n-2} \nabla u\ \cdot \nabla \varphi\\
&\leq& C\ \brac{\Vert H \Vert_\infty\ \Vert \nabla \varphi \Vert_{n,\log - \beta} + \Vert \varphi \nabla H \Vert_{n,\log - \beta}}\ \Vert \nabla u \Vert_{n}^{n-1} \Vert \nabla u \Vert_{n,\log \alpha}. 
\end{ma}
\]
Then the same line of arguments as before implies regularity.\\
\\
As for the \emph{structure} of this note: In Section~\ref{sec:orlicz} we will repeat most of the necessary results for logarithmic and exponential Orlicz spaces, e.g., H\"older inequalities in Section~\ref{ss:upperboundyc}, the Maximal Theorem (and as a corollary the result in \cite[Theorem II.1]{CLMS93}) in Section~\ref{ss:maxth}, and Trudinger's inequality in Section~\ref{ss:expgrowth}. In Section~\ref{s:frehse} we show via Frehse's counterexample that our theorem is non-trivial. Finally, the details of the above sketched arguments can be found in Section~\ref{s:pde}.\\
\\
The \emph{notation} we use is quite standard. For $p \in [1,\infty]$ we write $p'$ for the H\"older-conjugate exponent $p' = \frac{p}{p-1}$. We denote constants by $C$ and these depend on the quantities indexed by a subscript. We make no effort whatsoever to pick optimal constants, and in particular, these constants may vary from line to line. As several of the appearing constants depend on the dimensions involved, we do not especially denote this by a subscript. Finally, for two quantities $A,B \in \R$ we say that $A \aleq B$, if there is a positive constant $C$ such that $A \leq C\ B$. We write $A \ageq B$ iff $B \aleq A$ and say $A \aeq B$ iff $A \aleq B$, $A \ageq B$.

\vspace{2ex}\noindent
{\bf Acknowledgement.} This note is based on research which was conducted while the author was visiting the University of Warsaw. The author likes to thank in particular (in alphabetical order) Pawe\l{} Goldstein, S\l{}awomir Kolasi\'{n}ski, Pawe\l{} Strzelecki, and Anna Zatorska-Goldstein for many interesting discussions and the warm welcome in frosty Poland.\\
Additionally, it is in order to thank Professor Heiko von der Mosel at RWTH Aachen University for the constant support for many years now. The stay in Poland was financed by Strzelecki's and von der Mosel's grants with the FNP and DFG, respectively.

%% file: orlicz.tex
\section{Preliminaries -- Some Facts about Orlicz-Spaces}
In this section, we state and -- for the convenience of the reader -- prove several results on Orlicz Spaces of the type $L^p\log^\sigma L$. We are confident, that all these results are known and possibly only special cases of much more general theorems, but we limited our attention to the setting we are interested in. We emphasize, however, that the Orlicz space which we denote by $L^p \log^\sigma L$ is in part of the literature known as $L^p \log^{p\sigma} L$. Generally, we follow the notation in \cite{IO07}.
\label{sec:orlicz}
\subsection{Definition of the Relevant Spaces}
Let for $\alpha \in \R$, $p \in [1,\infty)$
\[
 [f]_{p,\log \alpha, \Omega} \equiv [f]_{L^p \log^\alpha L (\Omega)} := \intl_\Omega \abs{f}^p \log^\alpha (e+\abs{f}),
\]
and
\[
 \Vert f \Vert_{p,\log \alpha, \Omega} \equiv \Vert f \Vert_{L^p \log^\alpha L (\Omega)} :=  \inf \left \{ \lambda > 0: \quad [\lambda^{-1} f]_{p,\log \alpha,\Omega} \leq 1 \right \}.
\]
Moreover, we set for $\sigma > 0$
\[
 [f]_{\EXP,\sigma,\Omega} := \intl_{\Omega} \brac{e^{\abs{f}^\sigma} - 1 },
\]
and
\[
 \Vert f \Vert_{\EXP,\sigma, \Omega} :=  \inf \left \{ \lambda > 0: \quad [\lambda^{-1} f]_{\EXP, \sigma ,\Omega} \leq 1 \right \}.
\]
One checks, that these ''norms`` in fact satisfy each norm-condition but possibly the triangular inequality. In Lemma \ref{la:faketriangineq} we will prove that for $X = \EXP,\sigma$ or $X = p,\log \alpha$ we at least have for any $\Omega \subset \R^n$
\begin{equation}\label{eq:faketriineq}
 \Vert f + g\Vert_{X,\Omega} \leq C_X\ \brac{\Vert f \Vert_{X,\Omega} + \Vert g \Vert_{X,\Omega}}.
\end{equation}

\begin{lemma}\label{la:notlap2cond}
Let $p \in [1,\infty)$, $\sigma \in \R$. Then there exists a constant $\Lambda \equiv \Lambda_{p,\sigma} \geq 2$ such that for any $t > 0$ and any $\varphi(t) := t^p \log^\sigma (e+t)$
\[
 \varphi(\Lambda^{-1} t) \leq \frac{1}{2}\ \varphi(t).
\]
In particular, for any $L > 0$ we have for some constant $C_{p,\sigma,L} > 0$
\[
C_{p,\sigma,L}^{-1} \leq \frac{\inf \left \{ \lambda > 0: \quad [\lambda^{-1} f]_{p,\log \sigma,\Omega} \leq L \right \}}{\Vert f \Vert_{p,\log \sigma, \Omega}} \leq C_{p,\sigma,L}
\]
The same holds for $\varphi(t) := e^{t^\sigma}-1$, $\sigma > 0$.\\
\end{lemma}
\begin{proofL}{\ref{la:notlap2cond}}
Let $\varphi(t) = t^p \log^{\sigma} (e+t)$ where $\sigma \in \R$. We have for any $\Lambda \geq 3$,
\[
 \varphi(\Lambda^{-1} t) = \Lambda^{-p} \brac{\frac{\log(e+t)}{\log(e+\Lambda^{-1}t)}}^{-\sigma} \varphi(t).
\]
If $\sigma \geq 0$, obviously, $\log(e+\Lambda^{-1} t) \leq \log (e+t)$ implies
\[
\varphi (\Lambda^{-1} t) \leq \Lambda^{-p} \varphi(t),
\]
and the claim is proven. We are thus left with the case $\sigma < 0$. If $t \geq \Lambda^2$,
\[
\begin{ma}
 \log(e+t) &\leq& \log(\Lambda (e+\Lambda^{-1} t))\\
&=& \frac{1}{2}\log \Lambda^2 + \log(e+\Lambda^{-1} t) \leq \frac{1}{2} \log (e+t) + \log(e+\Lambda^{-1} t),\\
\end{ma}
\]
and consequently,
\[
 \log(e+t) \leq 2 \log(e+\Lambda^{-1} t).
\]
Thus, for any $t \geq \Lambda^2$ and any $\Lambda \geq 3$
\[
 \varphi(\Lambda^{-1} t) = \Lambda^{-p}\ 2^{-\sigma}\ \varphi(t).
\]
If on the other hand $t \in (0,\Lambda^2)$,
\[
 \varphi(\Lambda^{-1} t) \leq \Lambda^{-p} \log^{-\sigma}(e+t)\ \varphi(t) \leq  \Lambda^{-p} \log^{-\sigma}(2 \Lambda^2)\ \varphi(t).
\]
We conclude by choosing $\Lambda \geq 3$ such that
\[
 \max(\Lambda^{-p} 2^{-\sigma}, \Lambda^{-p} \log^{-\sigma}(2 \Lambda^2)) \leq \frac{1}{2}.
\]
Finally, if for some $\sigma > 0$ we set $\varphi(t) = e^{t^\sigma} - 1$, we have
\[
 \varphi(2^{-\fracm{\sigma}} t) = \sum_{k=1}^\infty 2^{-k} t^{\sigma k} \leq 2^{-1} \sum_{k=1}^\infty t^{\sigma k} = 2^{-1} \varphi(t).
\]
\end{proofL}

\begin{lemma}\label{la:faketriangineq}
Let for $p\in [1,\infty)$, $\alpha \in \R$, $\sigma > 0$ the function $\varphi: [0,\infty) \to [0,\infty)$ be defined as
\begin{equation}\label{eq:varphipeqtplogalpha}
 \varphi(t) = t^p\ \log^\alpha (e+t),
\end{equation}
or
\begin{equation}\label{eq:varphipeqetsigmam1}
 \varphi(t) = e^{t^\sigma} - 1.
\end{equation}
Then, for a constant $C = C_{\alpha}$ or $C = C_{\sigma}$, for all $s,t \geq 0$
\begin{equation}\label{eq:varphipseudoconvex}
 \varphi\brac{\frac{s+t}{2}} \leq C \brac{\varphi(t)+ \varphi(s)}. 
\end{equation}
Morerover, there exists a constant $C_\varphi > 0$ such that for any $0 \leq t_1 \leq t_2 < \infty$ we have
\begin{equation}\label{eq:varphipseudomon}
 \varphi(t_1) \leq C_\varphi\ \varphi(t_2).
\end{equation}
In particular, in view of Lemma \ref{la:notlap2cond}, \eqref{eq:faketriineq} holds.
\end{lemma}
\begin{proofL}{\ref{la:faketriangineq}}
It suffices to prove \eqref{eq:varphipseudoconvex} where $0 \leq s \leq t$. Then, \eqref{eq:varphipseudoconvex} is obvious, if $\varphi$ is monotone rising, as then
\[
 \varphi\brac{\frac{s+t}{2}}  \leq \varphi\brac{t} \leq \varphi(t) + \varphi(s).
\]
Thus the claim for the exponential $\varphi$ in \eqref{eq:varphipeqetsigmam1} and in view of
\[
\begin{ma}
 &&\brac{t^p\ \log^\alpha (e+t)}'\\
&=& t^{p-1}\ \log^\alpha(e+t) \brac{p+\alpha\ \log^{-1} (e+t)\ \frac{t}{e+t}}\\
&\geq& t^{p-1}\ \log^\alpha(e+t) \brac{p+\min(\alpha,0)}
\end{ma}
\]
also the claim for $\varphi$ as in \eqref{eq:varphipeqtplogalpha} for $\alpha \geq -p$ is obvious. So assume now that $\alpha < 0$. We have for $t \geq 4$, $0 \leq s \leq t$
\[
 \log\brac{e+\frac{s+t}{2}} \geq \log\brac{\frac{e+t}{2}} = \log(e+t)-\log 2 \geq \log(e+t)-\frac{1}{2} \log (e+t).
\]
If on the other hand $0 \leq s \leq t \leq 4$,
\[
 \log\brac{e+\frac{s+t}{2}} \geq 1 \geq \frac{\log(e+t)}{\log(e+4)}.
\]
Thus, taking
\[
 C_\alpha := \brac{\max (\log (e+4),2)}^{-\alpha},
\]
for any $0 \leq s \leq t < \infty$, and $\alpha < 0$, using the convexity of $()^p$,
\[
 \varphi\brac{\frac{s+t}{2}} \leq C_\alpha\ \varphi(t) \leq C_\alpha\ \brac{\varphi(t) + \varphi(s)}.
\]
As for \eqref{eq:varphipseudomon}, this is obvious, if $\varphi$ is monotone, i.e. if $\varphi$ is like \eqref{eq:varphipeqetsigmam1} or for $\alpha \geq -p$ \eqref{eq:varphipeqtplogalpha}. So let $\alpha \leq p$, then
\[
 \varphi'(t) = t^{p-1}\ \log^\alpha(e+t)\ \brac{\frac{t}{e+t} \log^{-1}(e+t)\alpha + p}.
\]
In order for $\varphi'(t)$ to be negative, $t$ has to be such that
\[
 \frac{t}{e+t} \log^{-1}(e+t)\alpha + p < 0,
\]
and there exists $0 < a < b < \infty$ such that this possibly happens only if $t \in (a,b)$. If $t_1, t_2 > b$, the claim is obvious, because $\varphi(t)$ is monotone rising, as it is in the case $t_1,t_2 < a$. Assume now that $t_1 \leq b$ and $t_2 \geq a$, then
\[
 \varphi(t_1) \leq \frac{\max_{(0,b)} \varphi}{\min_{(a,\infty)} \varphi}\ \varphi(t_2).
\]
As $\varphi(t) \neq 0$ for all $t \neq 0$ and $\lim_{t \to 0} \varphi = 0$, $\lim_{t \to \infty} = \infty$, the constant in the last inequality is finite.
\end{proofL}

\begin{lemma}\label{la:logfhs}
Let $f \in L^1(\R^n)$. Then, for any $p \geq 1$, $s > 0$, $\sigma \in \R$ we have for a constant depending on $s$, $p$ and $\sigma$,
\[
 [f^s]_{p,\log \sigma,\R^n} \aeq [f]_{sp,\log \sigma,\R^n}.
\]
In particular,
\[
 \Vert f^s \Vert_{p,\log \sigma,\R^n} \aeq \Vert f \Vert_{sp,\log \sigma,\R^n}^s.
\]
\end{lemma}
\begin{proofL}{\ref{la:logfhs}}
We have for any choice of $s > 0$ that
\[
 \log (e+t^s) \aeq \log (e+t).
\]
In fact, applying l'H\^{o}pital's rule, we have
\[
 \lim_{t \to \infty} \frac{\log(e+t^s)}{\log(e+t)} = s \lim_{t \to \infty} \frac{1 + t^{-1} e}{1 + t^{-s} e} = s.
\]
As moreover,
\[
 \lim_{t \to 0} \frac{\log(e+t^s)}{\log(e+t)} = 1,
\]
and for any $t \in (1,\infty)$
\[
 \frac{\log(e+t^s)}{\log(e+t)} > 0,
\]
we conclude that for any $t > 0$
\[
 0 < c_s \leq \frac{\log(e+t^s)}{\log(e+t)} \leq C_s < \infty.
\]
Inserting this in the definition of $[f]_{p,\log \sigma,\R^n}$ the lemma is proven.
\end{proofL}

\subsection{Absolute Continuity of the Norms}
\begin{lemma}\label{la:absolutecontinuity}
Assume that $f \in L^1(\R^n)$ and for $\alpha \in \R$, $p \in (1,\infty)$
\[
 [f]_{p,\log \alpha,\R^n} < \infty.
\]
Then there is for any $\varepsilon > 0$ a constant $\delta > 0$ such that the following holds:
\[
 \Vert f \Vert_{p,\log \alpha,A} \leq \varepsilon, \quad \mbox{for any $A \subset \R^n$ such that $\abs{A} \leq \delta$}.
\]
\end{lemma}
\begin{proofL}{\ref{la:absolutecontinuity}}
As
\[
 [f]_{p,\log \alpha,\R^n} = \left \Vert \abs{f}^p \log^\alpha (e+\abs{f}) \right \Vert_{1,\R^n} < \infty,
\]
there exists for any $\tilde{\varepsilon} > 0$ a $\tilde{\delta} > 0$ such that whenever $A \subset \R^n$ and $\abs{A} \leq \tilde{\delta}$ then
\[
 [f]_{p,\log\alpha,A} = \left \Vert \abs{f}^p \log^\alpha (e+\abs{f}) \right \Vert_{1,A} < \tilde{\varepsilon}.
\]
It is now enough to show that $\Vert f \Vert_{p,\log \alpha,A}$ tends to zero as $\tilde{\varepsilon}$ goes to zero: For any $\lambda \in (0,1)$ such that (for sufficiently small $\tilde{\varepsilon}$)
\begin{equation}\label{eq:lambdaplogpbetweeneps}
 \begin{cases}
  \frac{1}{2\tilde{\varepsilon}}  \leq \lambda^{-p} \brac{1+\log \brac{\lambda^{-1}}}^\alpha \leq \frac{1}{\tilde{\varepsilon}} \quad &\mbox{if $\alpha \geq 0$},\\
\lambda := \tilde{\varepsilon}^{\frac{1}{p}} \quad & \mbox{if $\alpha < 0$,}
 \end{cases}
\end{equation}
we have
\[
[\lambda^{-1} f]_{p,\log\alpha,A} = \lambda^{-p} \intl_{A} \abs{f}^p\ \log^\alpha (e+\abs{f})\ \brac{\frac{\log(e+\lambda^{-1}\abs{f})}{\log(e+\abs{f})}}^\alpha.
\]
If $\alpha < 0$, this implies
\[
 [\lambda^{-1} f]_{p,\log\alpha,A} \leq \lambda^{-p}\ [f]_{p,\log\alpha,A} \overset{\eqref{eq:lambdaplogpbetweeneps}}{\leq} \frac{\tilde{\varepsilon}}{\tilde{\varepsilon}} = 1.
\]
Next, for any $t \geq 0$, $\lambda < 1$ we have
\[
 \log(e+\lambda^{-1}t) \leq \log(\lambda^{-1}(e+t)) = \log(e+t) + \log(\lambda^{-1}),
\]
so
\[
 \frac{\log(e+\lambda^{-1}t)}{\log(e+ t)}\leq 1 + \log(\lambda^{-1}).
\]
Consequently if $\alpha \geq 0$,
\[
\begin{ma}
[\lambda^{-1} f]_{p,\log \alpha,A} &\leq& \lambda^{-p} (1+\log \brac{\lambda^{-1}})^\alpha\ \intl_{A} \abs{f}^p\ \log^\alpha (e+\abs{f})\\
&\overset{\eqref{eq:lambdaplogpbetweeneps}}{\leq}& \frac{1}{\tilde{\varepsilon}} \tilde{\varepsilon} = 1.
\end{ma}
\]
That is, whatever the sign of $\alpha$ may be, for the respective choice of $\lambda$,
\[
 \Vert f \Vert_{p,\log \alpha,A} \leq \lambda.
\]
As obviously $\lambda \xrightarrow{\tilde{\varepsilon} \to 0} 0$, we are done.
\end{proofL}

\subsection{Upper Bounds for Young's Complement}\label{ss:upperboundyc}
In this section, we calculate some examples for H\"older inequality in the setting of Orlicz-spaces.
\begin{lemma}[where $r \neq p$]\label{la:ubycrneqp}
Let $1 \leq r < p < \infty$ and $\sigma,\alpha \in \R$, $\alpha > -r$. Then there is a constant $C_{r,p,\sigma,\alpha} > 0$ such that for any $s,t \in (0,\infty)$
\begin{equation}\label{eq:ubystrestneqp}
 (st)^r \log^\alpha (e+st) - s^p \log^\sigma (e+s) \leq C_{r,p,\sigma,\alpha}\ t^q \log^\gamma (e+t),
\end{equation}
where $q \in (r,\infty)$ such that $\fracm{r} = \fracm{p} + \fracm{q}$, and
\[
\frac{\alpha}{r} = \frac{\sigma}{p} + \frac{\gamma}{q}.
\]
In particular, cf. Lemma \ref{la:notlap2cond}, the constant can be chosen such that for any $f,g \in L^1(\Omega)$, $\Omega \subset \R^n$
\[
\Vert f g \Vert_{r,\log \alpha,\Omega} \leq C_{r,p,\alpha,\sigma}\ \Vert f \Vert_{p,\log\sigma,\Omega} \ \Vert g \Vert_{q,\log \gamma,\Omega}.
\]
\end{lemma}
\begin{proofL}{\ref{la:ubycrneqp}}
In this general case we came up with no better approach than brute force calculation of the supremum of the left hand side of \eqref{eq:ubystrestneqp}: Fix $t > 0$. Write the left hand side of \eqref{eq:ubystrestneqp} as
\begin{equation}\label{eq:ubybtsdef}
\begin{ma}
 B_t(s) &:=& (st)^r \log^\alpha (e+st) - s^p \log^\sigma (e+s)\\
&=& s^r \brac{ t^r \log^\alpha (e+st) - s^{p-r} \log^\sigma (e+s)}.
\end{ma}
\end{equation}
One checks that as $0 < r < p$ -- whatever the sign of $\alpha$, $\sigma$ may be -- the limit $\lim_{s \to 0} B_t(s) = 0$ and $\lim_{s \to \infty} B_t(s) = -\infty$. Thus, if the supremum is not attained at $s = 0$, i.e. $B_t(s) = 0$, it is attained at some point $s \in (0,\infty)$ where
\[
\dot{B}_t(s) \equiv \frac{d}{ds}B_t(s) = 0.
\]
This is equivalent to
\[
\begin{ma}
&&r s^{r-1} t^r \log^\alpha (e+st) + \alpha \frac{s^r t^{r+1}}{e+st}\ \log^{\alpha - 1}(e+st) \\
&=& p s^{p-1} \log^\sigma(e+s) + \sigma \frac{s^p}{e+s}\ \log^{\sigma-1} (e+s).
\end{ma}
\]
Multiplying by $s^{1-r}$ (remember, $s > 0$) this is
\begin{equation}\label{eq:uybstep1}
\begin{ma}
&& t^r\ \log^\alpha (e+st) \brac{r  + \alpha \frac{s t}{e+st}\ \log^{-1}(e+st)}\\
&=&  s^{p-r} \log^\sigma(e+s) \brac{p+ \sigma \frac{s}{e+s}\ \log^{-1} (e+s)}.
\end{ma}
\end{equation}
One checks that
\[
 p+ \sigma \frac{s}{e+s}\ \log^{-1} (e+s) < \frac{1}{2}
\]
is equivalent to
\[
 \log (e+s) < \frac{-\sigma}{p-\frac{1}{2}} \frac{s}{e+s}.
\]
So whenever this happens, we know for constants $0 < c_{p,\sigma} < C_{p,\sigma}  < \infty$ that
\[
 s \in (c_{p,\sigma},C_{p,\sigma}).
\]
By the same reasoning, we can find $0 < c_{r,\alpha} < C_{r,\alpha} < \infty$ such that
\[
 r  + \alpha \frac{s t}{e+st}\ \log^{-1}(e+st) < \frac{1}{2}
\]
implies that
\begin{equation}\label{eq:ubystincc}
 st \in (c_{r,\alpha},C_{r,\alpha}).
\end{equation}
In particular, both sides of \eqref{eq:uybstep1} are positive but possibly for points
\[
t \in (c_{r,p,\alpha,\sigma}, C_{r,p,\alpha,\sigma}),\quad s \in (c_{p,\sigma}, C_{p,\sigma}),
\]
so the claim \eqref{eq:ubystrestneqp} holds for these $t,s$ by a convenient choice of the constant on the right-hand side. That is, we have to show the inequality \eqref{eq:ubystrestneqp} for points $s,t > 0$ such that (by \eqref{eq:uybstep1})
\begin{equation}\label{eq:ubyseqtlogst}
s =  t^{\frac{r}{p-r}}\ \log^{-\frac{\sigma}{p-r}}(e+s)\ \log^{\frac{\alpha}{p-r}} (e+st) \brac{\frac{r  + \alpha \frac{s t}{e+st}\ \log^{-1}(e+st)}{p+ \sigma \frac{s}{e+s}\ \log^{-1} (e+s)}}^{\frac{1}{p-r}},
\end{equation}
and we can assume that both numerator and denominator are positive. In particular, \eqref{eq:uybstep1} is then also equivalent to
\[
 \frac{r  + \alpha \frac{s t}{e+st}\ \log^{-1}(e+st)}{p+ \sigma \frac{s}{e+s}\ \log^{-1} (e+s)} = \frac{s^{p-r} \log^\sigma(e+s)}{t^r\ \log^\alpha (e+st)}.
\]
On the other hand, we assume that $B_{t}(s) > 0$ (if not, the claim is trivial), so in particular, cf. \eqref{eq:ubybtsdef},
\[
 t^r \log^\alpha (e+st) - s^{p-r} \log^\sigma (e+s) > 0.
\]
The last two facts imply that
\[
 \frac{r  + \alpha \frac{s t}{e+st}\ \log^{-1}(e+st)}{p+ \sigma \frac{s}{e+s}\ \log^{-1} (e+s)} \leq 1.
\]
Consequently, we can choose $\lambda = \lambda_{r,p,\alpha,\sigma} > 0$, such that
\[
s \overset{\eqref{eq:ubyseqtlogst}}{\leq}  t^{\frac{r}{p-r}}\ \log^{-\frac{\sigma}{p-r}}(e+s)\ \log^{\frac{\alpha}{p-r}} (e+st) \leq C_{r,p,\alpha,\sigma} \brac{1 + t^{\lambda} + s^{\frac{1}{2}}}.
\]
In particular,
\[
s \leq (e+t)^{C_\lambda},
\]
and
\begin{equation}\label{eq:ubylogepsleq}
 1 \leq \log (e+s),\ \log(e+st) \leq C\ \log(e+t).
\end{equation}
As $\alpha > -r$, we have that
\[
 r + \alpha \frac{st}{e+st}\log^{-1}(e+st) \geq r+\min(\alpha,0) > 0.
\]
Plugging this and \eqref{eq:ubylogepsleq} into \eqref{eq:ubyseqtlogst}, for our favorite $\varepsilon \in \brac{0,\frac{r}{p-r}}$ we get
\[
s \geq  t^{\frac{r}{p-r}}\ \log^{-\frac{\sigma}{p-r}}(e+s)\ \log^{\frac{\alpha}{p-r}} (e+st)\ c_{r,p,\alpha,\sigma} \overset{\eqref{eq:ubylogepsleq}}{\geq} c_{\varepsilon,r,p,\alpha,\sigma}\ t^{\frac{r}{p-r}-\varepsilon},
\]
so for constants depending on $r,p,\alpha,\sigma$,
\[
 \log (e+s), \log(e+st) \aeq  \log(e+t)
\]
Altogether, this implies
\[\begin{ma}
 B_{t} (s) &\leq& t^{\frac{r^2}{p-r}} \log^{-\sigma \frac{r}{p-r}}(e+s)\ \log^{\alpha \frac{r}{p-r}}(e+st)\ t^r \log^\alpha(e+st)\\
&\aeq & t^q\ \log^{-\sigma \frac{r}{p-r} + \alpha \frac{r}{p-r}+\alpha}(e+t).
\end{ma}
\]
\end{proofL}
The following might not be optimal...
\begin{lemma}[where $r = p$]\label{la:dualexplnlogalv2}
Let $p \in [1,\infty)$, $\alpha \in [0,p)$, $\sigma > 0$. Then for a constant $C_{p,\alpha,\sigma} > 0$ and any $s,t > 0$ we have
\[
 (st)^p \log^{-\alpha} (e+st) - s^p \log^\sigma(e+s) \leq C_{p,\alpha,\sigma}\ \chi_{t > 1}\ \exp \brac{C_{p,\alpha,\sigma}\ t^{\gamma}},
\]
where
\[
 \gamma = \frac{p}{\sigma+\alpha}.
\]
In particular, the constant can be chosen such that for any $f,g \in L^1(\Omega)$
\[
\Vert f g \Vert_{p,\log(-\alpha),\Omega} \leq C_{r,p,\alpha,\sigma}\ \Vert f \Vert_{p,\log\sigma,\Omega} \ \Vert g \Vert_{\EXP,\gamma,\Omega}.
\]
\end{lemma}
\begin{proofL}{\ref{la:dualexplnlogalv2}}
Again, brute force: The claim is obvious for $t \leq 1$, as then for any $s \geq 0$
\[
 t^p \log^{-\alpha} (e+st) - \log^\sigma(e+s) \leq 0.
\]
It thus suffices to prove the claim for $t > \bar{t}$, where $\bar{t}$ is chosen depending only on $p,\alpha,\sigma$. First of all, let $\bar{t} > 1$ to be chosen later. Then fix $t > \bar{t}$ and let
\[
 B_t(s) := (st)^p \log^{-\alpha} (e+st) - s^p \log^\sigma(e+s).
\]
We know that the maximum is attained at a point $s \in (0,\infty)$, as for $s = 0$ we have $B_t(0) = 0$ and for $s \to \infty$, $B_t(\infty) = - \infty$. Then our $s$ is such that $\frac{d}{ds} B_t(s) = 0$, i.e.
\[
\begin{ma}
&& p s^{p-1} t^p \log^{-\alpha}(e+st) - \alpha s^{p-1} t^p \log^{-\alpha - 1} (e+st) \frac{st}{e+st}\\
&=& p s^{p-1} \log^\sigma (e+s) + \sigma s^{p-1} \log^{\sigma-1}(e+s) \frac{s}{e+s},
\end{ma}
\]
which is equivalent to
\[
\begin{ma}
 &&t^p \log^{-\alpha} (e+st) \brac{p - \alpha \frac{st}{e+st}\ \log^{-1}(e+st)}\\
&=& \log^\sigma(e+s) \brac{p + \sigma \frac{s}{e+s} \log^{-1}(e+s)}.
\end{ma}
\]
As $\sigma > 0$ and $s >0$, both sides have to be  strictly positive, and we can rewrite the equation as
\begin{equation}\label{eq:expyoungquotient}
\frac{t^p \log^{-\alpha} (e+st)}{\log^\sigma(e+s)}  =  \frac{p + \sigma \frac{s}{e+s} \log^{-1}(e+s)}{p - \alpha \frac{st}{e+st}\ \log^{-1}(e+st)}.
\end{equation}
As $s$ is supposed to be a point of maximal value of $B_t$ and thus $B_t(s) \geq 0$, we know that
\[
 t^p \log^{-\alpha} (e+st) \geq \log^{\sigma}(e+s),
\]
and thus
\[
 \frac{p + \sigma \frac{s}{e+s} \log^{-1}(e+s)}{p - \alpha \frac{st}{e+st}\ \log^{-1}(e+st)} \geq 1,
\]
and of course we have an estimate in the other direction as well,
\[
 \frac{p + \sigma \frac{s}{e+s} \log^{-1}(e+s)}{p - \alpha \frac{st}{e+st}\ \log^{-1}(e+st)} \leq \frac{p+\sigma}{p-\alpha} =: \Lambda.
\]
Consequently, by \eqref{eq:expyoungquotient},
\[
 \frac{t^p \log^{-\alpha} (e+st)}{\log^\sigma(e+s)} \in (1,\Lambda).
\]
This justifies the following estimates: Firstly,
\begin{equation}\label{eq:youngexpstep1}
 e+s \leq \exp\brac{t^{\frac{p}{\sigma}} \log^{-\frac{\alpha}{\sigma}} (e+st)}.
\end{equation}
Secondly, we have also
\[
 \log (e+st) \geq \Lambda^{-\frac{1}{\sigma}}\ t^{\frac{p}{\sigma}} \log^{-\frac{\alpha}{\sigma}} (e+st),
\]
i.e.
\[
 \log^{\frac{\sigma+\alpha}{\sigma}} (e+st) \geq \Lambda^{-\frac{1}{\sigma}}\ t^{\frac{p}{\sigma}},
\]
which implies
\begin{equation}\label{eq:youngexpstep2}
 \log^{-\frac{\alpha}{\sigma}} (e+st) \leq \Lambda^{\frac{\alpha}{\sigma(\sigma+\alpha)}}\ t^{-\frac{p}{\sigma}\frac{\alpha}{\sigma+\alpha}}.
\end{equation}
Thus,
\[
 \begin{ma}
s &\overset{\eqref{eq:youngexpstep1}}{\leq}& \exp\brac{t^{\frac{p}{\sigma}} \log^{-\frac{\alpha}{\sigma}} (e+st)}\\
&\overset{\eqref{eq:youngexpstep2}}{\leq}&  \exp\brac{\Lambda^{\frac{\alpha}{\sigma(\sigma+\alpha)}}\ t^{\frac{p}{\sigma}\brac{1-\frac{\alpha}{\sigma+\alpha}}}}.
 \end{ma}
\]
Now,
\[
 \gamma := \frac{p}{\sigma}\brac{1-\frac{\alpha}{\sigma+\alpha}}> 0.
\]
Finally, we can conclude 
\[
 B_t(s) \leq s^p t^p \leq \exp\brac{C_{p,\sigma,\alpha} t^\gamma}\ t^p \overset{t > \bar{t}}{\leq} \exp\brac{C_{p,\sigma,\alpha}\ t^\gamma}.
\]
\end{proofL}

\subsection{Maximal Theorem for logarithmic Orlicz-Spaces}\label{ss:maxth}
We introduce the Hardy-Littlewood maximal function for measurable $f: \R^n \to \R$
\[
 \mathcal{M}f(x) := \sup_{r > 0}\ \abs{B_r}^{-1} \intl_{B_r(x)} \abs{f(y)}\ dy, \qquad \mbox{for $x \in \R^n$}.
\]
The following result is then well known:
\begin{lemma}[Maximal Theorem]\label{la:maxthm} (see, e.g., \cite[Theorem I.3.1]{Stein93}0)\\
There is a constant $C > 0$ such that for any $f \in L^1(\R^n)$ and any $s > 0$
\[
 s\ \abs{\{\abs{\mathcal{M}f} > s\}} \leq C\ \Vert f \Vert_{1}.
\]
\end{lemma}
The following corollary is a special case of the results obtained in \cite{KT82}, \cite{BK94}.
\begin{corollary}[Maximal Theorem for logarithmic Orlicz-spaces]\label{co:maxllogl}
Let $p \in (1,\infty)$, $\alpha \in (-p,\infty)$. Then for a constant $C_{p,\alpha}$ and for all $f \in L^1(\R^n)$,
\[
 [\mathcal{M} f]_{p,\log \alpha,\R^n} \leq C_{p,\alpha} [f]_{p,\log \alpha,\R^n}.
\]
In particular, in view of Lemma \ref{la:notlap2cond}, one can choose the constant so that also
\[
 \Vert \mathcal{M} f \Vert_{p,\log\alpha,\R^n} \leq C_{p,\alpha} \Vert f \Vert_{p,\log\alpha,\R^n}.
\]
\end{corollary}
\begin{proofC}{\ref{co:maxllogl}}
For arbitrary $s > 0$ we decompose
\[
 f_s(x) := f(x) \chi_{\{\abs{f} > s\}} .
\]
With this notation one checks
\[
\begin{ma}
 \abs{\{\mathcal{M} f > 2s \}} &\leq& \abs{\{\mathcal{M} f_s > s \}} + \abs{\{\mathcal{M} (f-f_s) > s \}}\\
&=& \abs{\{\mathcal{M} f_s > s \}}\\
&\overset{\sref{L}{la:maxthm}}{\leq}& C\ s^{-1}\ \Vert f_s \Vert_{L^1}\\
&=& C\ s^{-1}\intl_{\R^n} \chi_{\{\abs{f} > s\}} \abs{f}.
\end{ma}
\]
Thus, setting $\varphi(t) := t^p \log^\alpha (e+t)$, using that $\varphi' \geq 0$ as $\alpha \geq -p$
\[
 \begin{ma}
 \intl_{\R^n} \varphi (\mathcal{M} f(x))\ dx &\overset{\varphi(0)=0}{=}& \intl_{0}^\infty \varphi'(s) \abs{\{\mathcal{M} f > s\}}\ ds\\
&\overset{\varphi' \geq 0}{\leq}&  C\intl_{0}^\infty \varphi'(s)\ s^{-1} \intl_{\R^n} \chi_{\{\abs{f} > \frac{s}{2}\}} \abs{f(x)}\ dx\ ds\\
&=&  C \intl_{\R^n} \abs{f(x)} \intl_{0}^{2\abs{f}} \varphi'(s) s^{-1}\ ds\ dx\\
&\overset{p>1}{=}&  C \intl_{\R^n} \frac{1}{2} \varphi (2\abs{f}) + C \intl_{\R^n} \abs{f(x)} \intl_{0}^{2\abs{f}} \varphi(s) s^{-2}\ ds\ dx\\
 \end{ma}
\]
Now we claim, that our choice of $\alpha$, $p$ implies for any $t \geq 0$
\begin{equation}\label{eq:logexpdualintvpsm2gvptm1}
 \intl_{0}^{t} \varphi(s) s^{-2}\ ds \leq C_{\alpha,p}\ \varphi(t)\ t^{-1}.
\end{equation}
Observe, that because of $p > 1$,
\[
 s^{-1} \varphi(s) \xrightarrow{s \to \infty} \infty,\quad s^{-1 -\frac{p-1}{2}} \varphi(s) \xrightarrow{s \to 0} 0,
\]
and consequently the left hand side of \eqref{eq:logexpdualintvpsm2gvptm1} is strictly positive and bounded for any $t \in (0,\infty)$, tends to zero as $t$ tends to zero and blows up as $t$ tends to infinity. The same asymptotics hold for the right-hand side of \eqref{eq:logexpdualintvpsm2gvptm1}. Then, in order to show \eqref{eq:logexpdualintvpsm2gvptm1}, it suffices to use l'H\^{o}pital's rule and calculate the behavior at the extremal points $0$ and $\infty$:
\[
\lim_{t \to 0,\infty} \brac{\frac{\intl_{0}^{t} \varphi(s) s^{-2}\ ds}{\varphi(t) t^{-1}}} = \lim_{t \to 0,\infty} \frac{\varphi(t) t^{-2}}{\brac{\varphi(t) t^{-1}}'}.
\]
Now we have
\[
 \begin{ma}
\brac{\varphi(t) t^{-1}}' &=& \brac{t^{p-1} \log^\alpha (e+t)}'\\
&=& t^{p-2} \log^\alpha (e+t) \brac{p-1 + \alpha \frac{t}{e+t} \log^{-1} (e+t)}.
\end{ma}
\]
That is,
\[
 \brac{\varphi(t) t^{-1}}' = \frac{\varphi(t)}{t^2}\ \brac{p-1 + \alpha \frac{t}{e+t} \log^{-1} (e+t)}.
\]
As $p > 1$, we can assure that $\brac{p-1 + \alpha \frac{t}{e+t} \log^{-1} (e+t)} > 0$ for all $t$ sufficiently close to $0$ or $\infty$. Thus,
\[
 \lim_{t \to 0,\infty} \brac{\frac{\intl_{0}^{t} \varphi(s) s^{-2}\ ds}{\varphi(t) t^{-1}}} = \frac{1}{p-1}.
\]
On the other hand, this quotient is a smooth function on $(0,\infty)$, so we can conclude that the quotient is bounded by a constant depending on $p$ and $\alpha$.
\end{proofC}

Once we have the Maximal Theorem, we have of course an extension of the famous result by Coifman, Lions, Meyer, and Semmes relating the Hardy space $\mathcal{H}$ and div-curl-terms, see \cite[Theorem II.1]{CLMS93}.
\begin{theorem}\label{th:clms}
Let $E,F \in L^1(\R^n,\R^n)$, $n \geq 2$, and assume
\[
 \dv (E) = 0, \qquad \curl(F) = 0 \quad \mbox{weakly on $\R^n$.}
\]
Then for any $p > 1$ and
\[
-\min\brac{p,1+\frac{p}{n}} < \alpha < \min\brac{p,p-1+\frac{p}{n}}
\]
the following estimate holds:
\[
 \Vert E\cdot F \Vert_{\mathcal{H}(\R^n)} \leq C_{p,\alpha}\ \Vert E \Vert_{p,\log \alpha,\R^n}\ \Vert F \Vert_{p',\log -\frac{\alpha}{p-1},\R^n}.
\]
\end{theorem}
\begin{proofT}{\ref{th:clms}}
By the estimates in \cite[Lemma II.1]{CLMS93} we have for any $s \in (1,p)$, $t \in (1,p')$ such that $\frac{1}{s} + \frac{1}{t} = 1 + \frac{1}{n}$,
\[
 \Vert E\cdot F \Vert_{\mathcal{H}(\R^n)}  \leq C_{s,t}\ \left \Vert \brac{\mathcal{M}(\abs{E}^s)}^{\frac{1}{s}}\ \brac{\mathcal{M}(\abs{F}^t)}^{\frac{1}{t}}\right \Vert_{1,\R^n}.
\]
By Lemma \ref{la:ubycrneqp} we have
\[
\begin{ma}
 &&\left \Vert \brac{\mathcal{M}(\abs{E}^s)}^{\frac{1}{s}}\ \brac{\mathcal{M}(\abs{F}^t)}^{\frac{1}{t}} \right \Vert_{1,\R^n}\\
&\leq& C_{p,\alpha}\ \left \Vert \brac{\mathcal{M}(\abs{E}^s)}^{\frac{1}{s}} \right \Vert_{p,\log\alpha,\R^n}\ \left \Vert \brac{\mathcal{M}(\abs{F}^t)}^{\frac{1}{t}} \right \Vert_{p',\log-\frac{\alpha}{p-1},\R^n}\\
&\overset{\sref{L}{la:logfhs}}{\aeq }& C_{p,\alpha} \left \Vert \mathcal{M}(\abs{E}^s) \right \Vert_{\frac{p}{s},\log\alpha, \R^n}^{\fracm{s}}\ \left \Vert \mathcal{M}(\abs{F}^t) \right \Vert_{\frac{p'}{\fracm{t}},\log-\frac{\alpha}{p-1},\R^n}^t\\
\end{ma}
\]
Assume first that $0 \leq \alpha < \min(p,p-1+\frac{p}{n})$, in particular $\frac{\alpha}{p-1} < 1 + \frac{p'}{n}$ and $\frac{\alpha}{p-1} < p'$. So the three intervals
\[
 \brac{\frac{\alpha}{p-1},1+\frac{p'}{n}} \cap \brac{1,p'} \cap \brac{\frac{p'}{n},p'} \neq \emptyset.
\]
Consequently, it is possible to choose $t \in (1,p')$, $t < n$, so that
\[
 \frac{\alpha}{p-1} < \frac{p'}{t} < 1 + \frac{p'}{n}.
\]
Observe,
\[
1 > \frac{1}{s} = 1 +\frac{1}{n} - \frac{1}{t} > \frac{1}{p}.
\]
Then by Corollary \ref{co:maxllogl}, on the one hand
\[
 \Vert \mathcal{M}(\abs{F}^t)\Vert_{\frac{p'}{t},\log-\frac{\alpha}{p-1},\R^n}^{\frac{1}{t}} \aleq  \Vert \abs{F}^t\Vert_{\frac{p'}{t},\log-\frac{\alpha}{p-1},\R^n}^{\frac{1}{t}} \overset{\sref{L}{la:logfhs}}{\aeq }  \Vert F\Vert_{p',\log-\frac{\alpha}{p-1},\R^n}.
\]
And on the other hand,
\[
 \Vert \mathcal{M}(\abs{E}^s) \Vert_{\frac{p}{s},\log\alpha,\R^n}^s \aleq  \Vert E \Vert_{p,\log\alpha,\R^n}.
\]
The same holds in the case $-\min(p,1+\frac{p}{n}) < \alpha < 0$: we can choose $s \in (1,p)$, $s < n$, so that
\[
 \alpha > - \frac{p}{s} > - 1 - \frac{p}{n},
\]
and $t \in (1,p')$, $t < n$.
\end{proofT}

\subsection{Poincar\'e Inequality for logarithmic Orlicz-spaces}
For $r > 0$, $x \in \R^n$, $w : \R^n \to \R$ let
\[
 (w)_{r,x} := \abs{B_r(x)}^{-1}\ \intl_{B_r(x)} w.
\]
As the scaling arguments for Orlicz-spaces seem to be rather unpleasant, for the here necessary Poincare's inequality, we will follow the strategy in \cite{Bjoern10}.
\begin{lemma}[Poincar\'e-Inequality]  \label{la:orliczpoinc}
For any $\alpha \geq 0$, $p \in [1,\infty)$, there is a constant $C_\alpha$ such that whenever $B_r(x) \subset \R^n$, $r \in (0,\infty)$, and $w \in W^{1,1}(\R^n)$,
\[
 \left [ r^{-1}\ \brac{w-(w)_{r,x}}\right ]_{p,\log\alpha,B_r(x)} \leq C_{p,\alpha}\ [C\ \nabla w]_{p,\log \alpha,B_{5r}(x)}
\]
that is
\[
 \Vert w-(w)_{r,x} \Vert_{p,\log \alpha,B_r(x)} \leq C_{\alpha}\ r\ \Vert \nabla w \Vert_{n,\log \alpha,B_r(x)}.
\]
\end{lemma}
\begin{proofL}{\ref{la:orliczpoinc}}
First of all, for any $y \in B_r(x)$, we have
\[
 \abs{w(z)-(w)_{r,x}} \aleq{} r^{1-n}\ \intl_{B_r(x)}\ \intl_{0}^1\ \abs{\nabla w(z+t (\xi-z))}\ dt\ d\xi.
\]
Using transformation rule $\zeta := z + t (\xi -z)$ and observing that $\abs{\zeta - z} = t \abs{\xi - z} \leq 2tr$,
\[
\begin{ma}
 \abs{w(z)-(w)_{r,x}} &\aleq{}& t^{-n} r^{1-n}\ \intl_{0}^1\ \intl_{B_{2tr}(z)} \abs{\nabla w(\zeta)}\ d\zeta\ dt\\
&=& t^{-n} r^{1-n}\ \intl_{0}^1\ \intl_{B_{2tr}(z)} \abs{\nabla w(\zeta)}\ \chi_{B_{3r}(x)}(\zeta)\ d\zeta\ dt\\
&\aleq{}& r^{1}\ \mathcal{M} \brac{\abs{\nabla w}\ \chi_{B_{3r}(x)}} (z).
\end{ma}
\]
In particular,
\[
 r^{-1}\ \brac{w-(w)_{r,x}} \aleq{} \mathcal{M} \brac{\abs{\nabla w}\ \chi_{B_{3r}(x)}} (z).
\]
Next, let $\varphi(t)$ be defined as
\[
 \varphi(t) := t^p\ \log^\alpha(e+t), \quad \mbox{$t \in [0,\infty)$}
\]
Note that $\varphi$ is monotone (we don't even need convexity), as $\alpha \geq 0$. In fact,
\[
 \varphi'(t) = t^{p-1} \log^\alpha (e+t)\ \brac{p+\alpha \frac{t}{e+t} \log^{-1} (e+t)} \geq 0,
\]
Consequently,
\[
  \left [ r^{-1}\ \brac{w-(w)_{r,x}}\right ]_{p,\log\alpha,B_r(x)} \leq  \left [ \mathcal{M} \brac{C \abs{\nabla w}\ \chi_{B_{3r}(x)}} \right ]_{p,\log\alpha,\R^n}.
\]
The Maximal Theorem, Corollary \ref{co:maxllogl}, allows us to estimate further
\[
\begin{ma}
 &&\left [ r^{-1}\ \brac{w-(w)_{r,x}}\right ]_{p,\log\alpha,B_r(x)}\\
&\leq& C_{p,\alpha} \left [ C \abs{\nabla w}\ \chi_{B_{3r}(x)} \right ]_{p,\log\alpha,\R^n}\\
&=& C_{p,\alpha} \left [ C\ \abs{\nabla w} \right ]_{p,\log\alpha,B_{3r}(x)}\\
\end{ma}
\]
\end{proofL}

\subsection{Trudinger's inequality}\label{ss:expgrowth}
In this section, we will repeat the famous estimate between a space of type $\EXP$ and $W^{1,n}_0$. A more general result of the same type, involving logarithmic Orlicz-spaces, can be found in \cite{FLS96}. The following is a consequence of the Sobolev-Poincar\'e-inequalities.
\begin{lemma}(see \cite[Theorem 7.15]{GT83})\\
Let $\Omega \subset \R^n$ and $u \in W^{1,n}_0(\Omega)$, then for uniform constants $c_1,c_2 > 0$
\[
 \left [ \frac{u}{c_1\ \Vert \nabla u \Vert_{n,\Omega}}\right ]_{\EXP,\frac{n}{n-1},\Omega} \leq c_2\ \abs{\Omega}
\]
\end{lemma}
Thus, in view of Lemma \ref{la:notlap2cond}, we obtain
\begin{theorem}\label{th:expgegennablau}
For $\Omega \subset \R^n$ there exists a constant $C_{\abs{\Omega}} > 0$, such that for any $u \in W^{1,n}_0(\Omega)$ we have
\[
 \Vert u \Vert_{\EXP,\frac{n}{n-1},\Omega} \leq C_{\abs{\Omega}}\ \Vert \nabla u \Vert_{n,\Omega}.
\]
In particular for any $R > 0$ there is a constant $C_R > 0$ such that, if $\Omega = B_r(x)$, $r \in (0,R)$, $x \in \R^n$, we have for any $u \in W^{1,n}_0(B_r(x))$
\[
\Vert u \Vert_{\EXP,\frac{n}{n-1},B_r(x)} \leq C_R\ \Vert \nabla u \Vert_{n,B_r(x)}. 
\]
\end{theorem}

%% file: frehse.tex
\section{Frehse's Counterexample}\label{s:frehse}
As usual in these contexts, we use the famous example by Frehse, \cite{Frehse73}, which shows that our supposed integrability condition for the solution $u$ does not rule out singularities a priori, that is, the structure of the PDE has crucial influence on the regularity of its solution $u$. Let
\[
 u(x) := \log \log \frac{4}{\abs{x}}, \qquad \mbox{$x \in B_1(0)$}.
\]
It is not difficult to show that for any $\varepsilon > 0$
\[
 \nabla u \in L^{n} \log^{n-1-\varepsilon} L (B_1(0)).
\]
In fact,
\[
 \abs{\nabla u(x)} = \brac{\abs{x} \log\frac{4}{\abs{x}}}^{-1},
\]
which implies
\[
\begin{ma}
 &&\intl_{B_1(0)} \abs{\nabla u(x)}^n \log^{n-1-\varepsilon} (e+\abs{\nabla u(x)})\ dx\\
&\aeq & \intl_{r = 0}^1 r^{-1} \log^{-n}\brac{\frac{4}{r}}\ \log^{n-1-\varepsilon} \brac{e+\frac{1}{r \log \frac{4}{r}}}\ dr.
\end{ma}
\]
As (by l'H\^{o}pital's rule)
\[
\begin{ma}
 &&\lim_{r \to 0} \frac{\log \brac{e+\frac{1}{r \log \frac{4}{r}}}}{\log \brac{e+\frac{1}{r }}}\\
&=& \lim_{r \to 0} \frac{e r + 1}{e r \log \frac{4}{r} + 1}\ \brac{1 - \frac{1}{\log\brac{\frac{4}{r}}}}\\
&=& 1,
\end{ma}
\]
we have for any $r \in [0,1]$
\[
 \log \brac{e+\frac{1}{r \log \frac{4}{r}}} \aeq  \log \brac{e+\frac{1}{r }}.
\]
By the same argument, for any $r \in [0,1]$
\[
 \log \frac{4}{r} \aeq  \log \brac{e+\frac{1}{r}}.
\]
Thus,
\[
\begin{ma}
 &&\intl_{B_1(0)} \abs{\nabla u(x)}^n \log^{n-1-\varepsilon} (e+\abs{\nabla u(x)})\ dx\\
&\aeq & \intl_{r = 0}^1 r^{-1} \log^{-n}\brac{\frac{4}{r}}\ \log^{n-1-\varepsilon} \brac{\frac{4}{r}}\ dr\\
&\overset{t = \log\brac{\frac{4}{r}}}{=}& \intl_{t = \log 4}^\infty t^{-1-\varepsilon} dt\\
&=& \frac{\brac{\log 4}^{-\varepsilon}}{\varepsilon}.
\end{ma}
\]

%% file: pdes.tex
\section{The PDE Estimates}\label{s:pde}
In this section, we give some more details of the argument leading to regularity, as sketched in the introduction.
\subsection{The Estimates}
\begin{lemma}[Local control of Gradient]\label{la:loccongrad}
There is a uniform constant $\tau \in (0,1)$ such that the following holds: For any $u \in W^{1,n}(\R^n)$ and any ball $B_r \equiv B_r(x) \subset \R^n$ and $a \in W^{1,n'}_0(B_{2r}(x))$ such that
\begin{equation}\label{eq:pdest:lapaeqdvnu}
 \begin{cases}
 \lap a = \dn{\nabla u} \quad &\mbox{in $B_{2r}(x)$},\\
  a = 0 \quad &\mbox{on $\partial B_{2r}(x)$},
 \end{cases}
\end{equation}
we have the following estimate (recall $\frac{1}{n} + \frac{1}{n'} = 1$)
\[
 \Vert \nabla u \Vert_{n,B_r}^n \leq \tau \Vert \nabla u \Vert_{n,B_{2r}}^n + \Vert \nabla a \Vert_{n',B_{2r}}^{n'}.
\]
\end{lemma}
\begin{proofL}{\ref{la:loccongrad}}
Let $\eta \in C_0^\infty(B_{2r})$ be the usual cutoff-function, that is $\eta \equiv 1$ on $B_{r}$, $\abs{\nabla \eta} \leq 2r^{-1}$. We set
\[
 v := \eta (u-(u)_r),
\]
where (in a slight abuse of prior notation) $(u)_r := \mvint_{B_{2r} \backslash B_r} u$. Then
\[
 \intl_{B_r} \abs{\nabla u}^n \leq \intl_{\R^n} \abs{\nabla v}^n
\]
and
\[
\begin{ma}
\intl_{\R^n} \abs{\nabla v}^n &=& \intl_{\R^n} \abs{\nabla u}^{n-2} \nabla u \cdot \nabla v + \intl_{\R^n} \brac{\abs{\nabla v}^{n-2} \nabla v -\abs{\nabla u}^{n-2} \nabla u} \cdot \nabla v\\
&=:& I + II.
\end{ma}
\]
As for $I$ we note that $v$ is an admissible testfunction for \eqref{eq:pdest:lapaeqdvnu}, and thus
\[
 I \leq \Vert \nabla a \Vert_{n',B_{2r}}\ \Vert \nabla v \Vert_{n,B_{2r}} \leq C_\varepsilon\ \Vert \nabla a \Vert_{n',B_{2r}}^{n'} + \varepsilon \Vert \nabla v \Vert_{n,B_{2r}}^n.
\]
Absorbing for small $\varepsilon > 0$, this implies
\[
 \intl_{\R^n} \abs{\nabla v}^n \leq C \Vert \nabla a \Vert_{n',B_{2r}}^{n'} + 2 \abs{II}.
\]
As for $II$, we have
\[
\begin{ma}
 &&\abs{\nabla v}^{n-2} \nabla v - \abs{\nabla u}^{n-2} \nabla u\\
&=& \abs{\nabla v}^{n-2} \nabla (v-u) + \brac{\abs{\nabla v}^{n-2} - \abs{\nabla u}^{n-2}} \nabla u.
\end{ma}
\]
As moreover
\[
 \abs{\nabla v}^{n-2} \leq \abs{\nabla (v - u)}^{n-2} + \abs{\nabla u}^{n-2} + \sum_{k=1}^{n-3} c_k \abs{\nabla u - \nabla v}^{k} \abs{\nabla u}^{n-2-k},
\]
we have
\[
 \brac{\abs{\nabla v}^{n-2} \nabla v - \abs{\nabla u}^{n-2}\nabla u}\cdot \nabla v  \leq C_\varepsilon\ \abs{\nabla (u-v)}^n + \varepsilon \brac{\abs{\nabla v}^{n} + \abs{\nabla u}^n}.
\]
Thus,
\[
 \abs{II} \leq \varepsilon \brac{\Vert \nabla v \Vert_{n,B_{2r}}^n + \Vert \nabla u \Vert_{n,B_{2r}}^n}+ C_\varepsilon \Vert \nabla (u-v) \Vert_{n,B_{2r}}^n.
\]
Again, absorbing for small $\varepsilon > 0$ this leaves us with
\[
 \Vert \nabla v \Vert_{n,B_{2r}}^n \leq C_\varepsilon\ \Vert \nabla (u-v) \Vert_{n,B_{2r}}^n + 2\varepsilon \Vert \nabla u \Vert_{n,B_{2r}}^n + C \Vert \nabla a \Vert_{n',B_{2r}}^{n'}.
\]
Now,
\[
 \nabla(u - v) = \nabla \brac{(1-\eta_r) (u - (u)_r)}
\]
and (note that the support of $\nabla (u -v)$ is a subset of $B_{2r} \backslash B_r$) Poincar\'e's inequality implies that
\[
 \Vert \nabla (u-v) \Vert_{n,B_{2r}}^n \aleq{} \Vert \nabla u \Vert_{n,B_{2r}\backslash B_r}^n.
\]
Thus, we have shown that
\[
 \Vert \nabla u \Vert_{n,B_r}^n \leq C_\varepsilon \Vert \nabla u \Vert_{n,B_{2r}\backslash B_r}^n + 2\varepsilon \Vert \nabla u \Vert_{n,B_{2r}}^n + C\ \Vert \nabla a \Vert_{n',B_{2r}}^{n'},
\]
which, using Widman's holefilling trick, implies
\[
 (C_\varepsilon+1)\Vert \nabla u \Vert_{n,B_r}^n \leq (C_\varepsilon + 2 \varepsilon)\ \Vert \nabla u \Vert_{n,B_{2r}}^n + C\ \Vert \nabla a \Vert_{n',B_{2r}}^{n'},
\]
which -- if we set $\tau := \frac{C_\varepsilon + 2\varepsilon}{C_\varepsilon+1} \in (0,1)$ for some $\varepsilon \in (0,\frac{1}{4})$ -- implies
\[
\Vert \nabla u \Vert_{n,B_r}^n \leq \tau \Vert \nabla u \Vert_{n,B_{2r}}^n + \Vert \nabla a \Vert_{n',B_{2r}}^{n'}.
\]
\end{proofL}

\begin{lemma}\label{la:intest}
For any $\alpha \in [0,\min(2,n-1))$, $\sigma = n-1-\alpha$ and any $R > 0$ there exists a constant $C_{r,\alpha}$ such that the following holds:\\
Let $\varphi \in C_0^\infty(B_r(x))$, for $r \in (0,R)$ and $x \in \R^n$, $w_1,\ldots,w_n \in W^{1,n}(\R^n)$, $H \in L^\infty \cap W^{1,n} (\R^n)$. Then
\[
\begin{ma}
 &&\intl_{\R^n} \det (\nabla w_1,\ldots,\nabla w_n)\ H\ \varphi\\
&\leq& C_{R,\alpha} \Vert \nabla w_1 \Vert_{n,\log\alpha,B_{6r}}\ \Vert \nabla w_2 \Vert_{n,B_2r}\ldots \Vert \nabla w_{n} \Vert_{n,B_{2r}}\\
&&\quad \cdot \Vert \nabla \varphi \Vert_{n,B_r} \brac{\Vert H \Vert_{\infty}  + \Vert \nabla H \Vert_{n,\log\sigma,\R^n}}.
\end{ma}
\]
\end{lemma}
\begin{proofL}{\ref{la:intest}}
Let $\eta_r \in C_0^\infty(B_{2r}(x))$ be the usual cutoff-function which equals one on $B_r(x)$. Again we denote $(w)_{2r}$ to be the mean value on $B_{2r}$. We define
\[
 \tilde{w}_i := \eta_r (w_i-(w_i)_{2r}), \quad \mbox{$1 \leq i \leq n$}
\]
Then, using partial integration, the duality of Hardyspace and BMO, and finally $W^{1,n}$-Poincar\'e inequality,
\[
\begin{ma}
 I &:=& \intl_{\R^n} \det (\nabla w_1,\ldots,\nabla w_n)\ H\ \varphi\\
&=&  \intl_{\R^n} \det (\nabla \tilde{w}_1,\ldots,\nabla \tilde{w}_n)\ H\ \varphi\\
&=&  \intl_{\R^n} \det (\nabla \tilde{w}_1,\ldots,\nabla \tilde{w}_{n-1},\nabla (H\varphi))\ \tilde{w}_n\\
&\aleq{}& \left \Vert \det (\nabla \tilde{w}_1,\ldots,\nabla \tilde{w}_{n-1},\nabla (H\varphi)) \Vert_{\mathcal{H},\R^n}\ \right \Vert \nabla \tilde{w}_n \Vert_{n,\R^n}\\
&\aleq{}& \left \Vert \det (\nabla \tilde{w}_1,\ldots,\nabla \tilde{w}_{n-1},\nabla (H\varphi)) \Vert_{\mathcal{H},\R^n}\ \right \Vert \nabla w_n \Vert_{n,B_{2r}}.
\end{ma}
\]
Now one can rewrite
\[
 \det (\nabla \tilde{w}_1,\ldots,\nabla \tilde{w}_{n-1},\nabla (H\varphi)) = L(\tilde{w}_1,\ldots,\tilde{w}_{n-1}) \cdot \nabla (H\varphi),
\]
where $L$ is a multilinear operator which is divergence free (cf. \cite{CLMS93}, or \cite{Kol09}). Then by Theorem \ref{th:clms} for any $\alpha \in [0,2)$
\[
\begin{ma}
 &&\Vert \det (\nabla \tilde{w}_1,\ldots,\nabla \tilde{w}_{n-1},\nabla (H\varphi)) \Vert_{\mathcal{H}}\\
&\overset{\sref{T}{th:clms}}{\aleq{}}& \Vert \abs{\nabla \tilde{w}_1} \ldots \abs{\nabla \tilde{w}_{n-1}}\Vert_{n',\log \frac{\alpha}{n-1},\R^n}\ \Vert \nabla (H\varphi)) \Vert_{n,\log-\alpha,\R^n}\\
&\overset{\sref{L}{la:ubycrneqp}}{\aleq{}}& \Vert \nabla \tilde{w}_1 \Vert_{n,\log \alpha,\R^n}\ \Vert \abs{\nabla \tilde{w}_{2}} \ldots \abs{\nabla \tilde{w}_{n-1}}\Vert_{\frac{n}{n-2},\R^n}\ \Vert \nabla (H\varphi)) \Vert_{n,\log-\alpha,\R^n}\\
&\leq& \Vert \nabla \tilde{w}_1 \Vert_{n,\log \alpha,\R^n}\ \Vert \nabla \tilde{w}_{2}\Vert_{n,\R^n} \ldots \Vert \nabla \tilde{w}_{n-1}\Vert_{n,\R^n}\ \Vert \nabla (H\varphi)) \Vert_{n,\log-\alpha,\R^n}\\
\end{ma}
\]
Using Lemma \ref{la:faketriangineq} and Poincar\'e inequality for Orlicz-spaces as in Lemma \ref{la:orliczpoinc},
\[
\begin{ma}
 &&\Vert \nabla \tilde{w}_1 \Vert_{n,\log \alpha,\R^n}\\
&\aleq{}& \Vert \nabla w_1 \Vert_{n,\log \alpha,B_{2r}} + r^{-1}\ \Vert \brac{w_1-(w_1)_{2r}} \Vert_{n,\log \alpha,B_{2r}}\\
&\aleq{}& \Vert \nabla w_1 \Vert_{n,\log \alpha,B_{6r}}.
\end{ma}
\]
This and usual $W^{1,n}$-Poincar\'e inequality imply
\[
\begin{ma}
I &\aleq{}& \Vert \nabla w_1 \Vert_{n,\log\alpha,B_{6 r}}\ \Vert \nabla w_2 \Vert_{n,B_{2r}}\ldots \Vert \nabla w_{n} \Vert_{n,B_{2r}} \cdot\\ 
&&\quad \cdot \brac{\Vert H \Vert_{\infty} \Vert \nabla \varphi \Vert_{n} + \Vert \varphi \nabla H \Vert_{n,\log-\alpha}}.
\end{ma}
\]
Now we apply Lemma \ref{la:dualexplnlogalv2} to get
\[
 \Vert \varphi \nabla H \Vert_{n,\log-\alpha,\R^n} \leq \Vert \nabla H \Vert_{n,\log \sigma,\R^n}\ \Vert \varphi \Vert_{\EXP,\gamma,\R^n}.
\]
In order to apply Trudinger's inequality, we like to have $n' \overset{!}{=} \gamma = \frac{n}{\sigma+\alpha}$, that is $\sigma + \alpha = n-1$. Then, by Theorem \ref{th:expgegennablau} (note that $\varphi \in C_0^\infty(B_r)$)
\[
 \Vert \varphi \Vert_{\EXP,\gamma,B_{r}} \leq C_{R,\alpha}\ \Vert \nabla \varphi \Vert_{n,B_r}
\]
\end{proofL}

\subsection{The Conclusion: Proof of Theorem \ref{th:th}}
Lemma \ref{la:loccongrad} and Lemma \ref{la:intest} imply that (if $\alpha$ and $\sigma$ are appropriately chosen) for the solution $u$ of \eqref{eq:hmodel} and any ball $B_{12 r}(x)$, $r \in (0,1)$ within the domain of the PDE,
\[
 \Vert \nabla u \Vert_{n,B_r}^n \leq \brac{\tau + C_{H,\alpha,\sigma} \Vert \nabla u \Vert_{n, \log\sigma, B_{12 r}}^{n'}}\ \Vert \nabla u \Vert_{n,B_{4r}}^{n}.
\]
So for any $r$ small enough, applying Lemma \ref{la:absolutecontinuity}, we have
\[
 \Vert \nabla u \Vert_{n,B_r}^n \leq \frac{\tau + 1}{2} \Vert \nabla u \Vert_{n,B_{2r}}^{n}.
\]
Note that $\frac{\tau+1}{2} < 1$, so iterating this one obtains $\beta > 0$ such that for all sufficiently small $r > 0$
\[
 \Vert \nabla u \Vert_{n,B_r} \leq r^\beta,
\]
which implies that $u$ is H\"older-continuous, by Dirichlet growth theorem.